\documentclass[pre,superscriptaddress,twocolumn]{revtex4-1}
\usepackage{amsmath}
\usepackage{amssymb}
\usepackage{amsfonts}
\usepackage{amsthm}
\usepackage{bbm}
\usepackage{wrapfig}
\usepackage{bm}
\usepackage[mathscr]{eucal}
\usepackage[margin=1in]{geometry}
\usepackage{graphicx}
\usepackage{color}
\usepackage[colorlinks=true, allcolors=blue]{hyperref}
\usepackage{algorithm}
\usepackage{algpseudocode}
\usepackage{enumitem}

\begin{document}
\title{Computing distances on Riemann surfaces}

\author{Huck Stepanyants}
\email[corresponding author: ]{stepanyants.j@northeastern.edu}
\affiliation{Department of Physics, Northeastern University, Boston, Massachusetts 02115, USA}
\affiliation{Network Science Institute, Northeastern University, Boston, Massachusetts 02115, USA}

\author{Alan Beardon}
\affiliation{Department of Pure Mathematics and Mathematical Statistics, University of Cambridge, Cambridge CB3 0WB, UK}

\author{Jeremy Paton}
\affiliation{Department of Physics, Northeastern University, Boston, Massachusetts 02115, USA}
\affiliation{Network Science Institute, Northeastern University, Boston, Massachusetts 02115, USA}

\author{Dmitri Krioukov}
\affiliation{Department of Physics, Northeastern University, Boston, Massachusetts 02115, USA}
\affiliation{Network Science Institute, Northeastern University, Boston, Massachusetts 02115, USA}
\affiliation{Department of Mathematics, Northeastern University, Boston, Massachusetts 02115, USA}
\affiliation{Department of Electrical and Computer Engineering, Northeastern University, Boston, Massachusetts 02115, USA}

\begin{abstract}

Riemann surfaces are among the simplest and most basic geometric objects. They appear as key players in many branches of physics, mathematics, and other sciences. Despite their widespread significance, how to compute distances between pairs of points on compact Riemann surfaces is surprisingly unknown, unless the surface is a sphere or a torus. This is because on higher-genus surfaces, the distance formula involves an infimum over infinitely many terms, so it cannot be evaluated in practice. Here we derive a computable distance formula for a broad class of Riemann surfaces. The formula reduces the infimum to a minimum over an explicit set consisting of finitely many terms. We also develop a distance computation algorithm, which cannot be expressed as a formula, but which is more computationally efficient on surfaces with high genuses. We illustrate both the formula and the algorithm in application to generalized Bolza surfaces, which are a particular class of highly symmetric compact Riemann surfaces of any genus greater than~1.
\end{abstract}

\maketitle

\section{Introduction} \label{Sintroduction}

Riemann surfaces~\cite{Gilligan2012} are one of the simplest geometric objects, with widespread 
applications in physics, mathematics, network science, and other areas. In physics, Riemann surfaces appear naturally 
in dynamical systems. The simplest example is billiards~\cite{Masur,Gutkin}, i.e., 
a freely moving ball confined to a polygon-shaped table, and bouncing off its sides. 
This system is equivalently formulated as a particle 
moving along a straight line on a compact Riemann surface. In particular, a particle moving on the Bolza 
surface, a compact Riemann surface of genus~${g=2}$, is one of the simplest and earliest chaotic systems 
studied in physics~\cite{Hadamard}. 

Riemann surfaces also appear in integrable systems~\cite{Babelon,Gieres}. 
The trajectory of such a system in its parameter space is effectively confined to a Riemann surface. 
Elsewhere in mathematical physics, Riemann surfaces are a part of string theory~\cite{Kim,Bonelli}, 
potential theory~\cite{Gustafsson,Chirka}, and approximation theory~\cite{Komlov,Aptekarev1,Aptekarev2}. 
In quantum physics, Riemann surfaces appear in the study of fractional quantum Hall states~\cite{Wen1}, 
spin liquids~\cite{Wen2}, and quantum gravity~\cite{Aldrovandi,Distler}. 

In network science, Riemann surfaces have been used as latent spaces in some models of networks. 
For example, networks embedded in hyperbolic spaces are one of the first models that reproduce a 
collection of the most basic properties common to many real-world networks~\cite{Krioukov1}. 
Yet our main motivation for this paper comes from~\cite{Krioukov2,Hoorn2023}, where the convergence of the 
Ollivier curvature~\cite{Ollivier} of random geometric graphs~\cite{Penrose2003,Dall} to the Ricci curvature of their 
underlying Riemannian manifolds was proven. Vertices in such graphs are random points in a Riemannian 
manifold, and pairs of vertices are connected by edges if the distance between the two vertices in the 
manifold is below a fixed threshold. To build a random geometric graph in a manifold, one thus needs to 
be able to compute distances between pairs of points in it. 

\begin{figure*}[t!]
    \begin{center}
        \includegraphics[width=\textwidth]{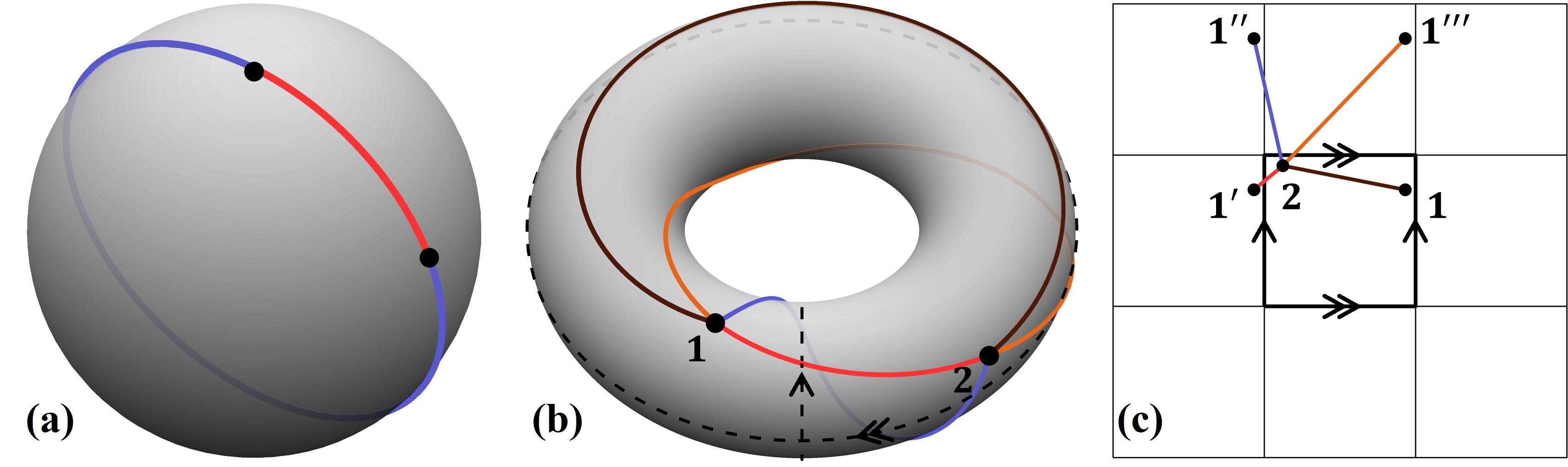}
        \caption{\textbf{(a)} The two geodesic paths connecting a given pair of points on a sphere. 
                \textbf{(b)} A topological representation of the 2D torus~${\mathbb{R}^2/\mathbb{Z}^2}$. 
                The dashed lines are its ``seams.'' Cutting the torus along these seams produces the 
                square in (c). The colored solid lines indicate four of the 
                infinitely many geodesic paths connecting the pair of points~${1,2}$. 
                \textbf{(c)} The corresponding points and geodesics in the~${\mathbb{R}^2}$ representation of the 
                torus. The highlighted central square is a fundamental domain of the torus and the arrows 
                indicate its side pairing. The distance between point 1 and 2 on the torus is between 2 and the 
                image~${1^\prime}$ of point~${1}$ on~${\mathbb{R}^2}$.}
        \label{F1}
    \end{center}
\end{figure*}

Compact Riemann surfaces are two-dimensional Riemannian manifolds. They can be equipped 
with a metric which defines the distance between a pair of points on the surface. 
Genus~${g=0}$ (sphere),~${g=1}$ (torus), and~${g \geq 2}$ surfaces admit the spherical, euclidean, and 
hyperbolic metrics, respectively. While computing distances on the first two of these surfaces 
is trivial, the task is much more difficult on a hyperbolic Riemann surface, 
where the distance is given as an infimum over a countably infinite number of geodesics, 
i.e., locally minimal paths. There are no formulas for the distance between points 
on any genus~${g \geq 2}$ surface which are computable in finite time. To be precise, while it is known 
that the mentioned infimum is actually a minimum since only a finite set of geodesics needs to be 
considered (Prop.~2.1,~\cite{Lu}), it is not known what this set actually is for 
any~${g \geq 2}$ surface, which is surprising since distance plays a crucial role in 
many applications of these surfaces mentioned above. 

Here we compute distances on Riemann surfaces obtained by pairing 
the sides of a convex hyperbolic polygon. We obtain an expression for a finite set of geodesics 
whose minimum gives the correct distance between any two points on these surfaces. This leads to a distance formula 
which can be evaluated in finite time. We also develop a fast algorithm for computing distances on these 
surfaces using a different set of ideas. For illustrative purposes, we show how the formula and the algorithm can be applied 
to computing distances on generalized Bolza surfaces~\cite{Ebbens}. These are highly symmetric hyperbolic Riemann surfaces 
of any genus~${g \geq 2}$. We prove that the time required to evaluate the formula in this case is~${O \left( g^{\alpha \log g} \right)}$, 
where~${\alpha \approx 1.52}$, and show experimentally that the time required to run the algorithm is~${O \left( g^4 \right)}$.

In Sec.~\ref{Sbackground}, we provide necessary background information on Riemann surfaces, and a detailed introduction to the distance problem. 
In Sec.~\ref{Sdistance_formula} and~\ref{Sdistance_algorithm}, we present and discuss our results, a formula and an algorithm, respectively, 
for computing distances on Riemann surfaces. In Sec.~\ref{SSg}, we apply our results to 
generalized Bolza surfaces, and evaluate the resulting running times. 

\section{Background and definitions} \label{Sbackground}

The distance between two points in a space is the length of the shortest geodesic, or locally minimal path, 
that connects them in the space. If the space is a compact Riemann surface, then the number of geodesics connecting a given pair of points depends 
on its \emph{genus}, which is the number of holes the surface has. Given two points on the genus~${0}$ sphere (Fig.~\ref{F1}(a)), there are 
always two geodesics between them, and thus the distance is simply a minimum of two lengths. In contrast, the number of geodesics between two 
points on the torus~${\mathbb{R}^2/\mathbb{Z}^2}$ is infinite (four are shown in Fig.~\ref{F1}(b)), 
which could make computing 
distances on the torus difficult. Nevertheless, due to its symmetry, the distance~${d}$ on the torus is given by the simple formula: 
\begin{align} \label{Edist_torus}
    \begin{split}
        d^2 = & \left( \frac{1}{2} - \Big\lvert \frac{1}{2} - |x_1-x_2| \Big\rvert \right)^2 \\
        + & \left( \frac{1}{2} - \Big\lvert \frac{1}{2} - |y_1-y_2| \Big\rvert \right)^2 ,
    \end{split}
\end{align}
where~${(x_1,y_1)}$ and~${(x_2,y_2) \in \left[ -1/2,+1/2 \right]^2}$ are the coordinates of the representatives of two points on the torus, which 
is the unit square with vertices at~${(\pm 1/2, \pm 1/2)}$ whose opposite sides are paired, Fig.~\ref{F1}(c). For genus~${g>1}$ (hyperbolic) surfaces, 
the number of geodesic paths connecting two points is also infinite. Yet there are no closed-form distance 
formulas analogous to Eq.~\eqref{Edist_torus} for any of these surfaces. 

Here, we address the distance problem for the hyperbolic Riemann surfaces known as \emph{quotient surfaces}. 
A quotient surface is constructed by taking the quotient of the hyperbolic plane~${\mathbb{H}^2}$ with respect to a discrete 
subgroup of its isometries, or distance-preserving transformations, called a \emph{Fuchsian group}, which we describe next. 

We will use the Poincar\'e disk model of~${\mathbb{H}^2}$, which is the 
complex unit disk equipped with the metric 
\begin{align}
    & ds^2 = \frac{4 dz d\bar{z}}{(1-|z|^2)^2}.
\end{align}
The group of orientation-preserving isometries of the Poincar\'e disk is the projective special 
unitary group~${\text{PSU}(1,1)}$, the quotient of the special unitary group~${\text{SU}(1,1)}$ by 
its center~${\{ I,-I \}}$. In matrix form, these isometries are 
\begin{equation} \label{E1}
    \gamma = \left[ \begin{matrix}
        a & \overline{c} \\ c & \overline{a}
    \end{matrix} \right], \quad
    \forall a,c \in \mathbb{C} ~:~
    \lvert a \rvert ^2 - \lvert c \rvert ^2 = 1 ,
\end{equation}
which act on points~${z \in \mathbb{C}}$ via fractional linear transformations, 
\begin{equation}
    \left[ \begin{matrix}
        a & \overline{c} \\ c & \overline{a}
    \end{matrix} \right]
    (z) =
    \frac{az + \overline{c}}{cz + \overline{a}} .
\end{equation}

\begin{figure}[t!]
    \begin{center}
        \includegraphics[width=0.5\textwidth]{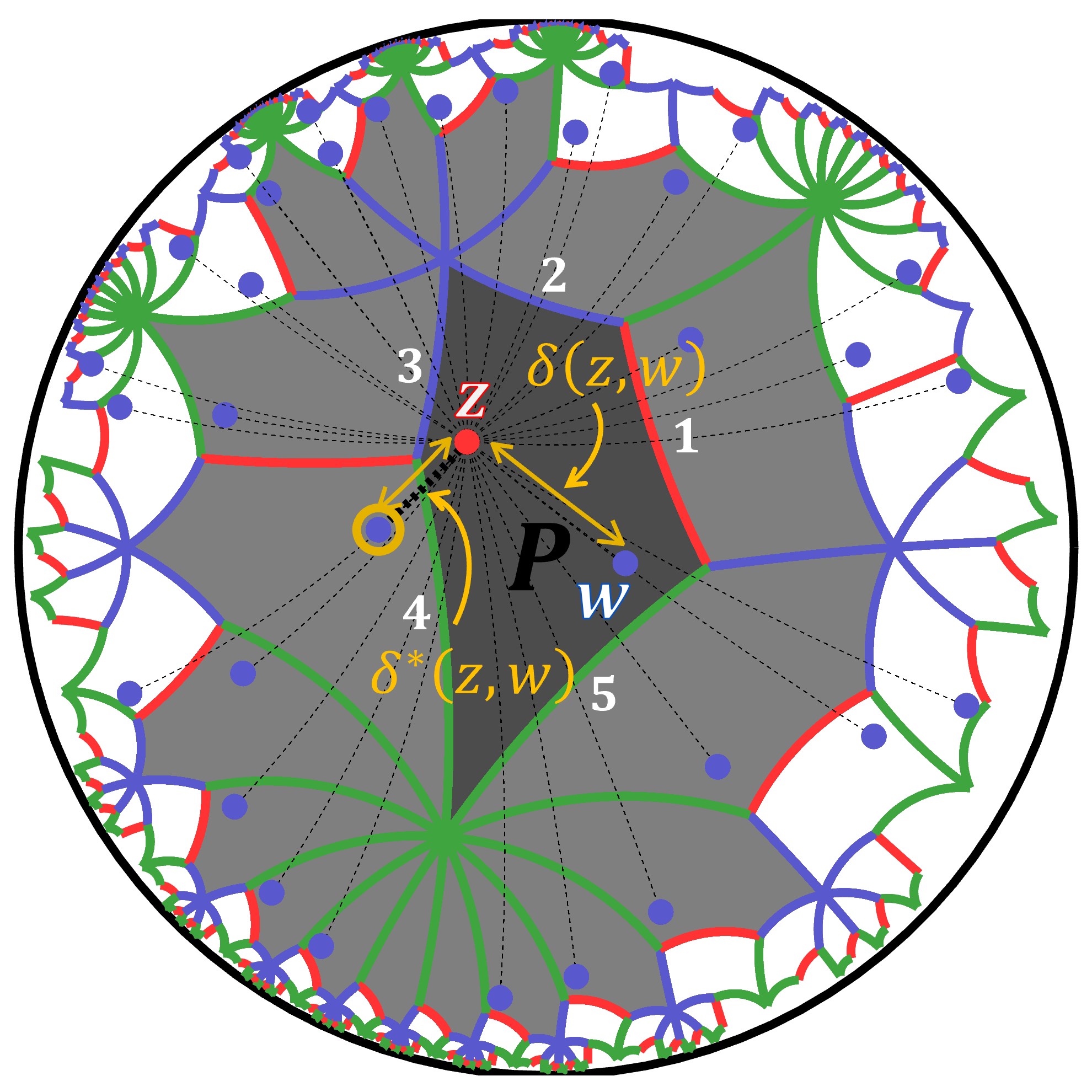}
        \caption{Quotient surface obtained by pairing the sides of a nonregular pentagon in the 
        order~${1 \longleftrightarrow 1 ,~ 2 \longleftrightarrow 3 ,~ 4 \longleftrightarrow 5}$. 
        The generators of the Fuchsian group~${\Gamma}$ map the fundamental polygon~${P}$ (dark 
        gray) to its edge-adjacent neighbors (light gray). To compute the 
        distance~${\delta^\star (z,w)}$ between a pair of points~${z}$ and~${w}$ on the surface, 
        the lengths~${\delta(z,\gamma w)}$ of the infinite number of geodesics between~${z}$ (red 
        dot) and all of the images~${\gamma w}$ (blue dots) of~${w}$,~${\forall \gamma \in \Gamma}$, 
        on the hyperbolic plane~${\mathbb{H}^2}$ must be considered. The shortest such distance, 
        between~${z}$ and the circled image of~${w}$ in the figure, is then~${\delta^\star(z,w)}$. }
        \label{F12}
    \end{center}
\end{figure}

A subgroup~${\Gamma}$ of the isometry group~${\text{PSU}(1,1)}$ 
is called \emph{Fuchsian} if it acts discontinuously 
on~${\mathbb{H}^2}$. This means that for any point~${z}$, the orbit~${\Gamma z}$ has no \emph{accumulation 
points}~${w ,~ |w| < 1}$. An accumulation point~${w}$ is one for which there are elements of~${\Gamma z}$ 
in any arbitrarily small punctured disk around~${w}$. Every Fuchsian group defines a quotient 
surface, 
\begin{align}
    & S = \mathbb{H}^2 \big\slash \Gamma .
\end{align}
Points~${[z]}$ on~${S}$ are orbits of points~${z}$ in~${\mathbb{H}^2}$ under actions by~${\Gamma}$, 
\begin{align}
    & [z] = \Gamma z .
\end{align}

The distance~${\delta^\star (z,w)}$ between two points~${[z]}$ and~${[w]}$ on~${S}$ is 
a distance between orbits~${\Gamma z}$ and~${\Gamma w}$ in~${\mathbb{H}^2}$. It is given by 
\begin{align} \label{Edstar}
    \begin{split}
        \delta^\star (z,w) & = \inf_{\gamma_1,\gamma_2 \in \Gamma} \delta(\gamma_1 z, \gamma_2 w) \\
        & = \inf_{\gamma \in \Gamma} \delta(z, \gamma w),
    \end{split}
\end{align}
where~${\delta}$ is the distance in~${\mathbb{H}^2}$. 
This formula cannot be directly evaluated, since~${\Gamma}$ has infinitely many 
elements. To find~${\delta^\star(z,w)}$, one must consider infinitely many geodesics from~${z}$ 
to~${\gamma w}$. Every such geodesic corresponds to a different minimal path between the 
same pair of points on~${S}$, like those illustrated in Fig.~\ref{F1} for the sphere and torus,
or in Fig.~\ref{F2} for the Bolza surface.

Here we will focus on surfaces~${S}$ whose Fuchsian groups~${\Gamma}$ are finitely generated and 
of the \emph{first kind}. The latter means that any point on the boundary of the unit disk, which is the 
boundary at infinity of the hyperbolic plane~${\mathbb{H}^2}$, is an accumulation point of some orbit 
of~${\Gamma}$. This guarantees that~${S}$ is a compact surface, and by 
Theorem~10.1.2 in~\cite{Beardon},~${S}$ 
has a convex \emph{fundamental} polygon~${P}$ with a finite number of sides. 
A fundamental polygon is one that contains 
exactly one representative of each point on~${S}$, see Fig.~\ref{F12}. It is known~\cite{Beardon} that 
in such settings, the generators of the group~${\Gamma}$ pair sides of~${P}$ by mapping it 
to its edge-adjacent neighbors, shown in light gray 
in Fig.~\ref{F12}. We will further assume that~${P}$ does not have any vertices at infinity. 

Using the geometry of~${P}$, we develop two methods for computing distances on~${S}$. 
In the next Sec.~\ref{Sdistance_formula}, we present a distance formula which reduces the 
infinities in Eq.~\eqref{Edstar} to finite sets. Then in Sec.~\ref{Sdistance_algorithm}, we develop 
a more efficient algorithmic approach to the problem. 

\begin{figure}[t!]
    \begin{center}
        \includegraphics[width=0.5\textwidth]{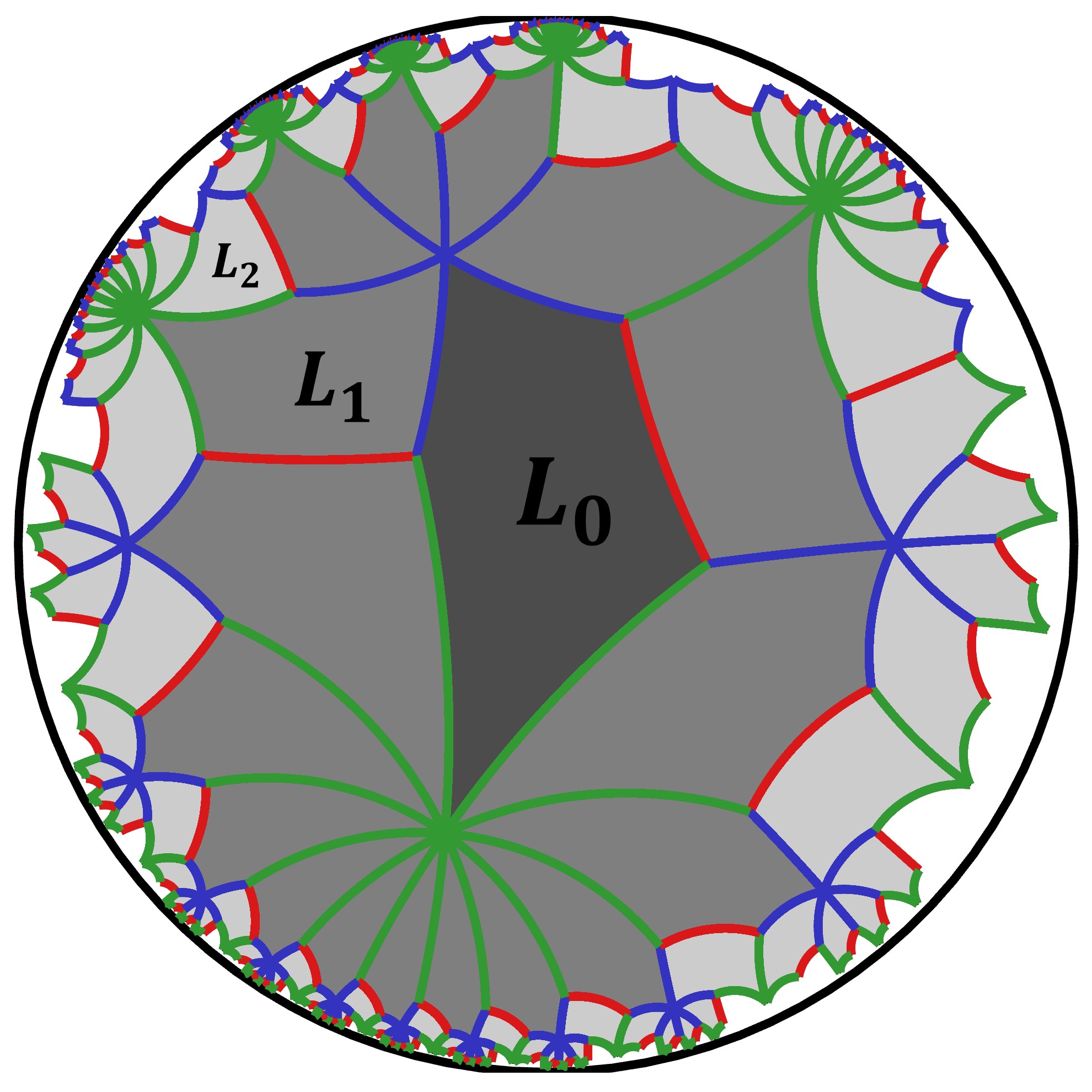}
        \caption{Shells~${L_0 = P}$ (dark gray),~${L_1}$ (medium gray), and~${L_2}$ (light gray). }
        \label{Fshells}
    \end{center}
\end{figure}

\section{Distance formula} \label{Sdistance_formula}

Let~${S}$ be a quotient surface obtained by pairing the sides of a convex fundamental 
polygon~${P}$, and let~${\Gamma}$ denote its Fuchsian group. Let~${z}$ and~${w}$ be 
representatives of two points on the surface. In this section we obtain a 
computable version of Eq.~\eqref{Edstar}, where instead of minimizing~${\delta(z, \gamma w)}$ 
over all~${\gamma \in \Gamma}$, we minimize over a finite subset~${\Gamma_0 \subset \Gamma}$. 
We will derive a formula for~${\Gamma_0}$ using the geometry of~${P}$. 

Let us call a pair of images of~${P}$ \emph{neighbors} if they share either a common edge or a common vertex. 
Let~${T \subset \Gamma}$ denote the subset of isometries from the Fuchsian group which map~${P}$ 
to all of its neighbors and itself. Then we define the~${k^\text{th}}$~\emph{shell}~${L_k}$, 
for~${k \geq 1}$, to be the difference 
\begin{align}
    & \bigcup_{\gamma \in T^k} \gamma P \bigg\backslash \bigcup_{\gamma \in T^{k-1}} \gamma P ,
\end{align}
and for~${k=0}$, we define~${L_0 = P}$, as illustrated in Fig.~\ref{Fshells}. 

\begin{figure}[t!]
    \begin{center}
        \includegraphics[width=0.5\textwidth]{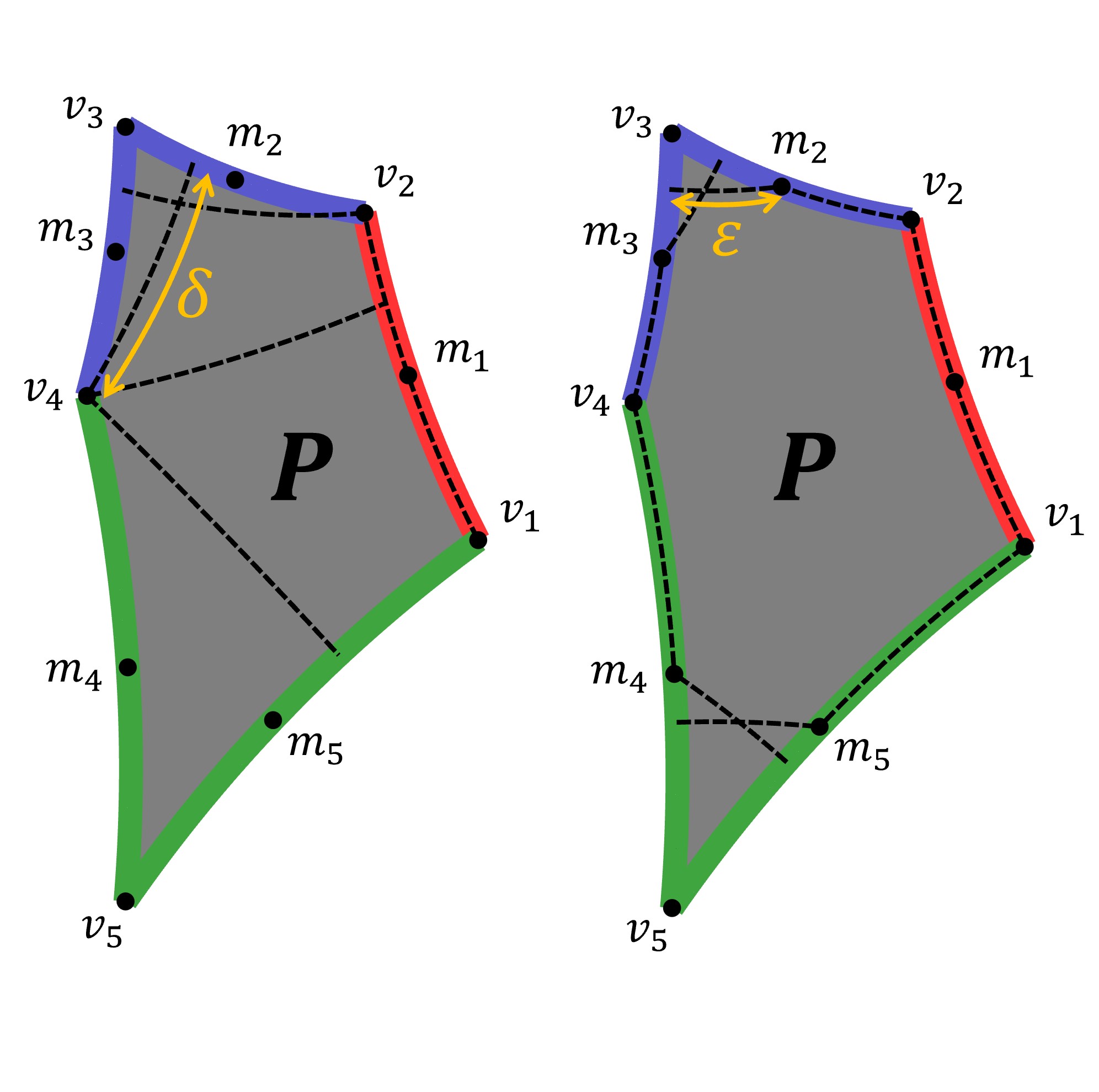}
        \caption{{\bf Left:} Shortest paths between nonadjacent sides whose minimum gives~${\delta}$. {\bf Right:} Shortest 
                paths from the midpoint of a side to an adjacent edge whose minimum gives~${\epsilon}$. }
        \label{Fepsilon_delta}
    \end{center}
\end{figure}

We next compute an upper bound on the minimum value of~${k}$,~${k_{\min}}$, for which the correct value 
of the distance~${\delta^\star(z,w)}$ is obtained by minimizing~${\delta(z, \gamma w)}$ 
over~${\forall \gamma \in T^k}$. Thus~${k_{\min}}$ is the minimum positive integer~${k}$ satisfying
\begin{align} \label{Ek}
    & \delta^\star(z,w) = \min_{\gamma \in T^k} \delta(z, \gamma w) 
\end{align}
for all~${z}$ and~${w}$. In other words, the finite subset~${\Gamma_0 \subset \Gamma}$ mentioned above 
can be taken to be 
\begin{align}
    & \Gamma_0 = T^{k_{\min}} .
\end{align}
We bound~${k_{\min}}$ for an arbitrary surface by the ratio of two distances: the surface 
diameter~${\mathscr{D}}$, i.e., the maximum possible distance between two points on the surface, and 
the minimal distance~${\mathscr{T}}$ required to traverse a shell. 
\emph{To traverse} here means \emph{to cross} from the outer boundary of~${L_k}$ to the outer 
boundary of~${L_{k+1}}$ for~${k \geq 0}$. In Appendix~\hyperref[AT_bound]{A}, we show that this distance does not depend 
on~${k}$, and satisfies the inequality 
\begin{align} \label{ET_upper_bound}
    \mathscr{T} \geq \min (\delta, 2\epsilon) ,
\end{align}
where~${\delta}$ is the minimum distance between any pair of non-adjacent sides of~${P}$, and~${\epsilon}$ is 
the minimum distance from the midpoint of an side to any adjacent side as illustrated in 
Fig.~\ref{Fepsilon_delta}. 

Then the distance required to cross from a starting point in~${L_0 = P}$ to a point 
in~${L_k}$ is lower-bounded by~${(k-1) \mathscr{T}}$. Since~${\mathscr{D}}$ is the surface diameter, 
then~${k_{\min}}$ must satisfy the inequality 
\begin{align} \label{Ekmin_intermediate}
    & (k_{\min} - 1) \mathscr{T} < \mathscr{D} .
\end{align}
But the surface diameter~${\mathscr{D}}$ is trivially upper bounded by the diameter 
of its fundamental polygon~${P}$,~${\text{diam} (P)}$: 
\begin{align} \label{ED_upper_bound}
    & \mathscr{D} \leq \text{diam} (P) = \max_{ij} \delta(v_i,v_j) ,
\end{align}
where the~${v_i}$ are the vertices of~${P}$. 
Combining Eqs.~\eqref{ET_upper_bound},~\eqref{Ekmin_intermediate}, and~\eqref{ED_upper_bound} leads 
to an upper bound on~${k_{\min}}$ in terms of geometric properties of~${P}$, which we denote~${k^\star}$: 
\begin{align} \label{Ekmin_bound}
    & k_{\min} \leq k^\star = \Bigg\lfloor 1 + \frac{\max_{ij} \delta(v_i,v_j)}{\min (\delta,2 \epsilon)} \Bigg\rfloor .
\end{align}
Substituting this bound in Eq.~\eqref{Ek} yields an explicit computable formula for the distance between 
a pair of points on the surface: 
\begin{align}
    & \delta^\star (z,w) = \min_{\gamma \in T^{k^\star}} \delta (z, \gamma w) .
\end{align}

\section{Distance algorithm} \label{Sdistance_algorithm}

In this section, we present an algorithm to compute~${\delta^\star(z,w)}$ 
by minimizing~${\delta(z,\gamma w)}$ over an even smaller subset of~${\Gamma}$. This subset depends on the location 
of~${z}$ and is unknown \emph{a priori}---instead, it is constructed via a search over the 
polygon tessellation. As we will see in Sec.~\ref{SSg}, this method is generally more computationally 
efficient than evaluating the distance formula in Sec.~\ref{Sdistance_formula}. 

\begin{figure}[t!]
    \begin{center}
        \includegraphics[width=0.45\textwidth]{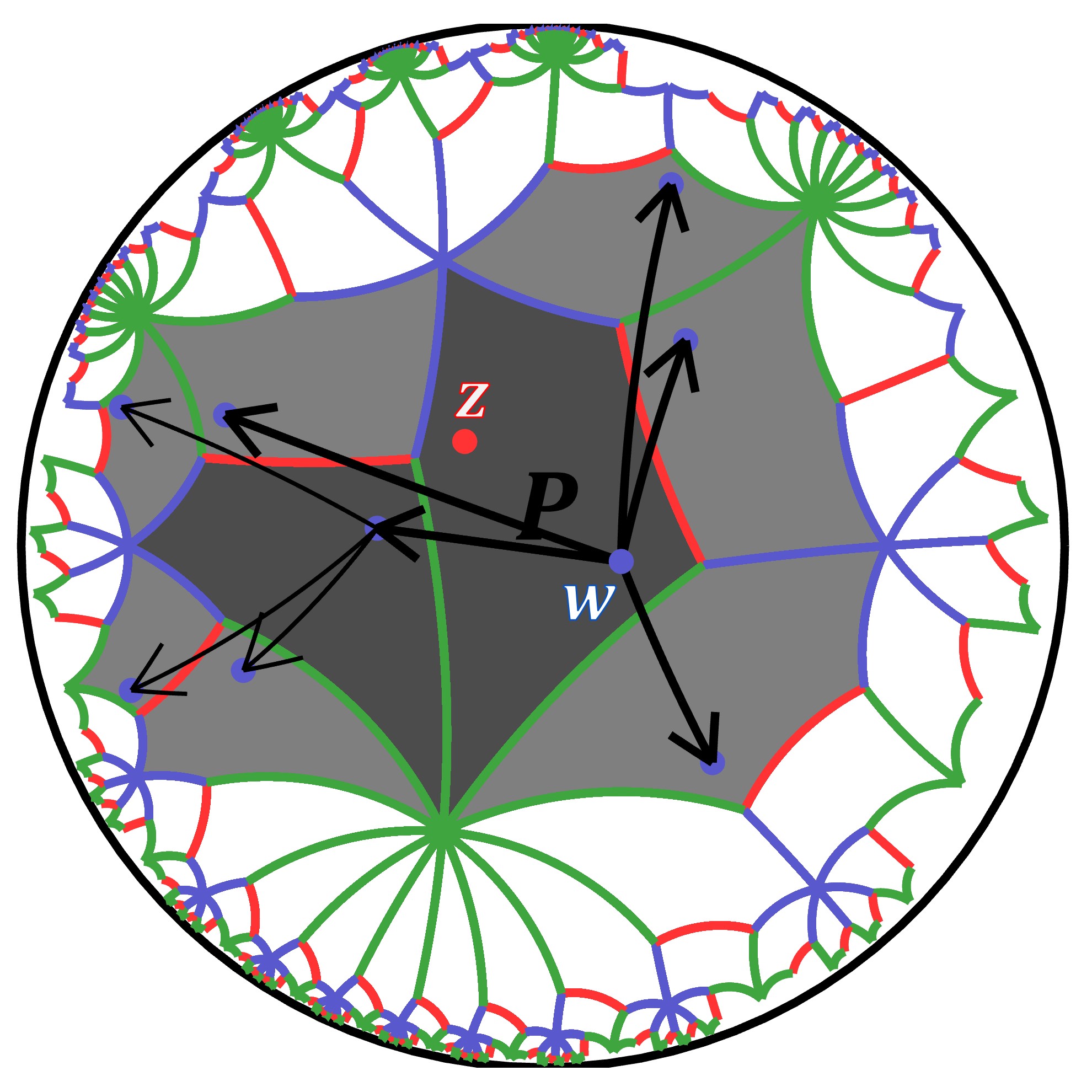}
        \caption{Illustration of one branch of the depth-first search algorithm. The search starts at 
        the fundamental polygon~${P}$ and progresses outward. At every step, the distance from~${z}$ to all 
        edge-adjacent polygons (light gray) of the last searched polygon (dark gray) are computed, 
        and those polygons are added to the list~${\Delta}$ if and only if they are at a greater 
        distance from~${z}$. The branch terminates when this distance exceeds the diameter of~$P$.}
        \label{Fsearch}
    \end{center}
\end{figure}

Let~${N}$ be the number of sides of~${P}$, and~${g_i,~ i = 1,2,\dots,N,}$ be the generators of~${\Gamma}$ 
and their inverses. With these notations, the algorithm is: 
\begin{algorithm}[H]
\caption{Distance calculation. }
\label{Adist}
\begin{algorithmic}
    \State${\Delta \gets \{1\}}$,~${\gamma_{\min} \gets 1}$
    \While{${|\Delta| > 0}$}
        \State${\gamma \gets \text{pop} (\Delta)}$
        \If{${\delta(z, \gamma w) < \delta(z, \gamma_{\min} w)}$}
            \State${\gamma_{\min} \gets \gamma}$
        \EndIf
        \For{${i = 1:N}$} 
            \If{${\delta(z, \gamma P) < \delta(z, \gamma g_i P) < \text{diam} (P)}$}
                \State${\Delta \gets \text{push}~(\gamma g_i)}$
            \EndIf
        \EndFor
    \EndWhile
    \State \textbf{return}~${\gamma_{\min}}$
\end{algorithmic}
\end{algorithm}
This algorithm is a depth-first search for an image~${\gamma w}$ of~${w}$, 
${\gamma \in \Gamma}$, that minimizes~${\delta(z,\gamma w)}$, in the polygon tessellation~${\Gamma P}$ 
illustrated in Fig.~\ref{Fsearch}. 
We keep track of two variables: the list of isometries~${\Delta \subset \Gamma}$ corresponding 
to the polygons to be searched next, and the current best solution~${\gamma_{\min}}$ 
that minimizes~${\delta(z,\gamma w)}$ across all~${\gamma}$ searched so far. 
Both~${\Delta}$ and~${\gamma_{\min}}$ are initialized to the 
identity:~${\Delta = \{1\}, \gamma_{\min} = 1}$, corresponding to~${P}$ itself. 

At every step of the search, we set~${\gamma}$ to be the last element 
of~${\Delta}$, and remove it from~${\Delta}$. 
If~${\delta(z, \gamma w) < \delta(z, \gamma_{\min} w)}$, we update~${\gamma_{\min}}$ 
to~${\gamma}$. We then append to~${\Delta}$ the group elements~${\gamma g_i}$,~${1 \leq i \leq N}$, 
corresponding to polygons that are edge-adjacent to the currently searched polygon, but only if they 
are at a greater distance from~${z}$: 
\begin{align} \label{Esearch_condition}
    & \delta(z, \gamma P) < \delta(z, \gamma g_i P) < \text{diam} (P) .
\end{align}
Here, the distances~${\delta(z,\gamma P)}$ and~${\delta(z,\gamma g_i P)}$ from~${z}$ 
to the polygons~${\gamma P}$ and~${\gamma g_i P}$ are defined to be the minimum distance 
from~${z}$ to one of their sides, with the exception~${\delta(z,P) = 0}$. 

The reason to consider only those polygons that increase the distance from~${z}$ is as follows. 
Consider the shortest path from~${z}$ to~${\gamma P}$, where~${\gamma \neq 1}$. It intersects an edge-adjacent 
polygon of~${\gamma P}$. This edge-adjacent polygon is closer to~${z}$ than~${\gamma P}$. Continuing in 
this way, we construct a chain of edge-adjacent polygons leading from~${\gamma P}$ to~${P}$ along which the 
distance to~${z}$ is decreasing, which guarantees that~${\gamma P}$ will be reached in the search. 

The reason not to search the polygons~${\gamma P}$ for 
which~${\delta (z, \gamma P)}$ exceeds~${\text{diam} (P)}$ is this: 
\begin{align}
    \begin{split}
        \delta (z, \gamma_{\min} P) & \leq \delta (z, \gamma_{\min} w) \\
        & \leq \mathscr{D} \\
        & \leq \text{diam} (P) ,
    \end{split}
\end{align}
where~${\gamma_{\min}}$ is the group element we are looking for, i.e., the one that 
minimizes~${\delta (z, \gamma w)}$, and~${\mathscr{D}}$ is the surface diameter. 

\begin{figure*}[t!]
    \begin{center}
        \includegraphics[width=1\textwidth]{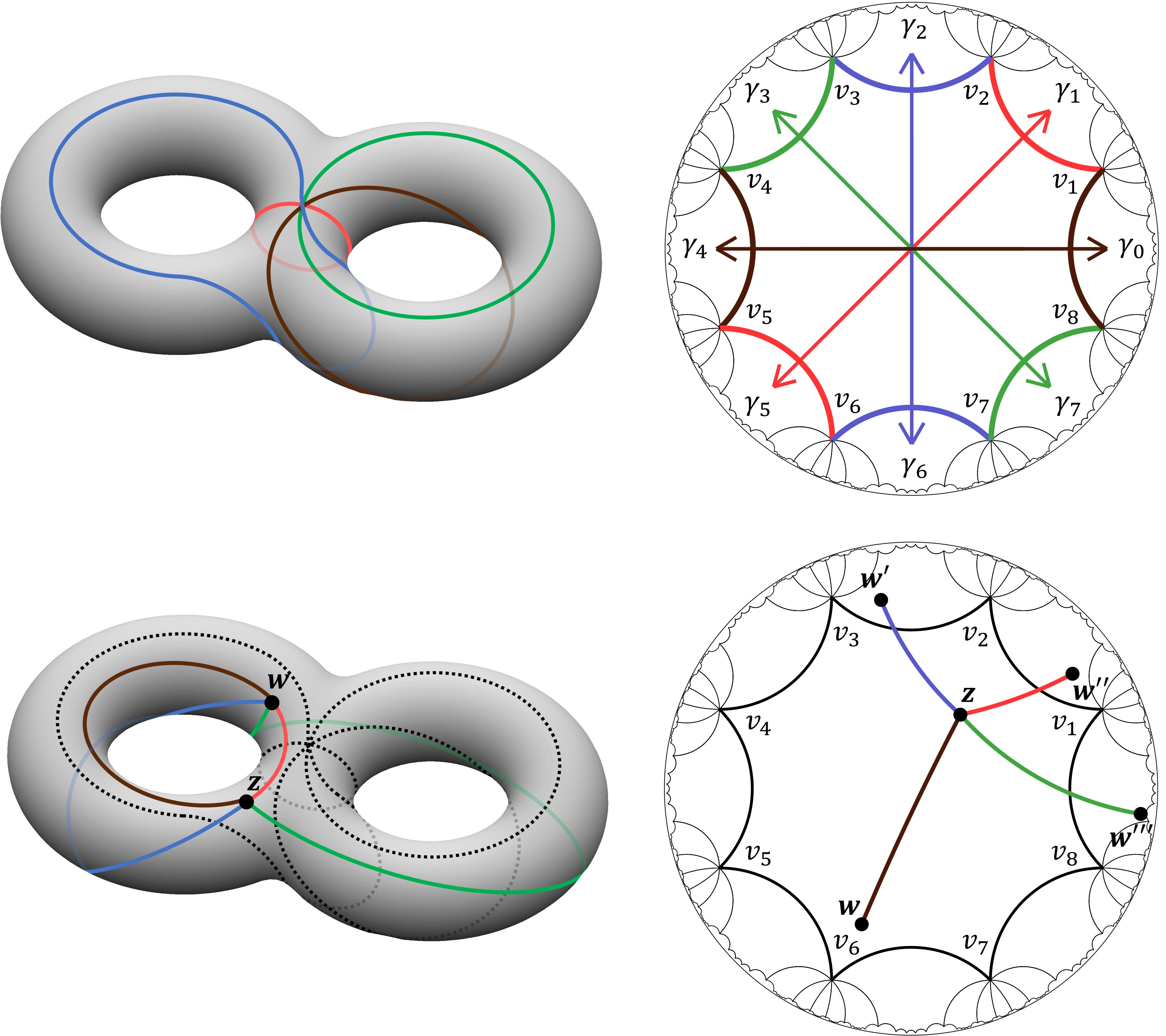}
        \caption{\textbf{Top Left:} A topological representation of the Bolza surface~${S_2}$. The colored lines indicate 
        its ``seams,'' cutting which and unfolding the surface produces the octagon on the right. 
        \textbf{Top Right:} The~${\mathbb{H}^2}$ representation of~${S_2}$. The colors indicate its side pairing, and 
        arrows show the action of the generators~${\gamma_i, i = 0,2,\dots,7}$ on the fundamental domain. 
        \textbf{Bottom Left:} Four of the infinitely many geodesic paths connecting two points on~${S_2}$. 
        The dashed lines are the seams at the top. 
        \textbf{Bottom Right:} The corresponding points and geodesics in the~${\mathbb{H}^2}$ representation of~${S_2}$. 
        The distance between points~${z}$ and~${w}$ on the surface is between~${z}$ and the 
        image~${w^{\prime\prime}}$ of~${w}$ on~${\mathbb{H}^2}$.}
        \label{F2}
    \end{center}
\end{figure*}

Finally, the condition in Eq.~\eqref{Esearch_condition} ensures that the search terminates 
after a finite time. To see this, let~${M}$ denote the number of polygons within the 
distance~${\text{diam} (P)}$ from~${z}$. Note that the number of branches of the search 
cannot exceed the number of sequences of such polygons ordered by increasing distance, which 
is~${2^{M}}$. Therefore, the number of polygons searched does not exceed~${M \cdot 2^{M}}$. 
We also consider a total of~${N}$ nearest neighbors of~${\gamma P}$ at every search step, even 
though not all of these neighbors are added to~${\Delta}$. Therefore the total number of steps in the 
algorithm is~${O(NM \cdot 2^M)}$. In practice, however, we find that the number of steps 
is~${O(NM)}$, implying that most polygons are searched only~${O(1)}$ times, as discussed in the next 
section. 

\section{Application to Generalized Bolza Surfaces} \label{SSg}

In this section, we apply our general results to the \emph{generalized Bolza surfaces}~\cite{Ebbens}, 
highly symmetric surfaces of any genus~${g \geq 2}$. We denote them by~${S_g}$. 
${S_2}$ is known as the Bolza surface~\cite{Bolza}. 

The generalized Bolza surface~${S_g}$ has as its fundamental polygon~${P}$ the regular~${4g}$-gon with 
interior angles~${\pi/2g}$ and vertices 
\begin{align}
    \begin{split}
        v_k & = \tanh (R/2) \exp \left( \left( k-\frac{1}{2} \right) \frac{\pi i}{2g} \right) ,\\
        k & = 0,1, \dots, 4g-1 ,
    \end{split}
\end{align}
where~${R}$ is the radius of~${P}$, which satisfies 
\begin{align}
    & \cosh R = \cot^2 \left( \frac{\pi}{4g} \right) .
\end{align}
The Fuchsian group~${\Gamma_g}$ of~${S_g}$ is generated by the following~${2g}$ isometries 
(Fig.~\ref{F2}), which pair opposite sides of the fundamental polygon: 
\begin{align}
    \begin{split}
        \gamma_k & = 
        \left[ \begin{matrix} 1 & \tanh \left( \frac{s}{2} \right) e^{\frac{k \pi i}{2g}} \\ 
        \tanh \left( \frac{s}{2} \right) e^{- \frac{k \pi i}{2g}} & 1\end{matrix} \right] ,\\
        k & = 0,1, \dots, 2g - 1 ,
    \end{split}
\end{align}
where~${s}$ is the length of a side of~${P}$, which satisfies 
\begin{align}
    & \cosh \left( \frac{s}{2} \right) = \cot \left( \frac{\pi}{4g} \right) .
\end{align}

Applying the distance formula to these surfaces, we observe that~${\delta < 2 \epsilon}$ 
(cf. Fig.~\ref{Fepsilon_delta}), and then upper bound~${\delta}$ 
using hyperbolic trigonometry to obtain 
\begin{align}
    & \min \{ \delta, 2 \epsilon \} = \delta < 
    \text{arccosh} \left( 2 \cos \left( \frac{\pi}{2g} \right) \right) .
\end{align}
The diameter~${\mathscr{D}_g}$ of~${S_g}$ is known~\cite{Stepanyants}, 
\begin{align}
    & \mathscr{D}_g = \text{arccosh} \left( \cot \left( \frac{\pi}{4g} \right) \right) ,
\end{align}
and therefore we can substitute the actual diameter, instead the upper bound in 
Eq.~\eqref{ED_upper_bound}, into Eq.~\eqref{Ekmin_bound} for~${k^\star}$ to obtain 
\begin{align}\label{Ekstar_Sg}
    & k^\star = \Bigg\lfloor 1 + \frac{\text{arccosh} \left( \cot \left( \frac{\pi}{4g} \right) \right)}{\text{arccosh} \left( 2 \cos \left( \frac{\pi}{2g} \right) \right)} \Bigg\rfloor .
\end{align}
The number of polygons in~${T^{k^\star}}$ is upper bounded by the number of sequences of~${2 k^\star}$ 
generators, 
\begin{align}
    & |T^{k^\star}| < (4g)^{2 k^\star} ,
\end{align}
and since~${k^\star \approx \log g / \text{arccosh}~2}$ for large~${g}$, the time required to evaluate the distance 
formula is~${O \left( g^{\alpha \log g} \right)}$, 
where~${\alpha = 2 / \text{arccosh}~2 \approx 1.52}$. 

For the Bolza surface~${S_2}$, evaluating Eq.~\eqref{Ekstar_Sg} yields~${k^\star = 2}$, but in Appendix~\hyperref[Sformula_Bolza]{B} 
we compute~${k_{\min}}$ exactly using different methods: 
\begin{align} \label{Ekmin_S2}
    & k_{\min} = 1. 
\end{align}

\noindent
Therefore, the number of~${P}$'s images one must search through to compute the distance 
is~${|T^{k_{\min}}| = 49}$, corresponding to all polygons adjacent to~${P}$ via either a 
side or a vertex. 

\begin{figure}[t!]
    \begin{center}
        \includegraphics[width=0.5\textwidth]{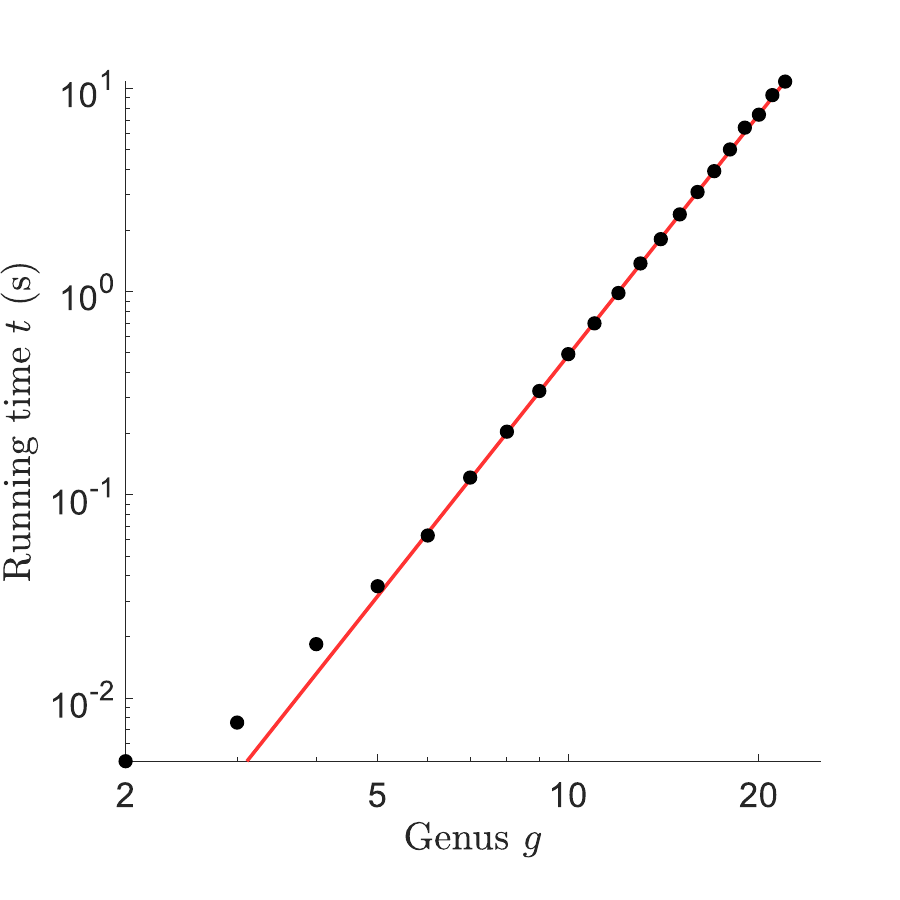}
        \caption{Distance algorithm running time~${t}$ on generalized Bolza surfaces~${S_g}$ 
        for~${g = 2,\dots,22}$. The red line shows the linear fit of the running times and 
        has slope~${3.94}$. }
        \label{FSg_running_time}
    \end{center}
\end{figure}

To apply the distance algorithm, we substitute the diameter of the surface for~${\text{diam} (P)}$ 
in Algorithm~\ref{Adist}, and find experimentally in Fig.~\ref{FSg_running_time} that the algorithm 
running time is~${O(g^4)}$. We cannot 
prove this observation because we cannot estimate the overcount---the average number of times the same 
polygon is searched. However, we can show that a maximum of~${O(g^4)}$ 
polygons can intersect a hyperbolic circle of radius~${R}$. Therefore, Fig.~\ref{FSg_running_time} 
suggests that each polygon is searched~${O(1)}$ times, so overcounting appears to be minimal. 

\appendix

\section{Lower bound on~${\mathscr{T}}$} \label{AT_bound}

In this section, we prove an upper bound on the shell-crossing distance~${\mathscr{T}}$ 
for any quotient surface~${S}$ obtained by pairing the sides of a convex polygon~${P}$. 

This distance is 
\begin{align}
    \nonumber & \mathscr{T} = \min \bigg\{ \delta \left( \partial T^n(P) , \partial T^{n+1}(P) \right) ,\quad n \geq 0 \bigg\}, 
\end{align}
where~${L_n = T^n (P)}$ is the~${n^\text{th}}$ shell,~${T \subseteq \Gamma}$ is the set consisting of all maps from~${P}$ to 
an edge- or vertex-adjacent polygon and the identity, and~${\partial L_n}$ is the boundary of~${L_n}$. 

First, observe that for~${\gamma \in T^n}$ we have
\begin{align}
    \nonumber & \gamma(P) \subseteq \gamma T (P) \subseteq L_{n+1} .
\end{align}
Hence, 
\begin{align}
    \nonumber & \delta \Bigl( \partial \gamma(P), \partial \left( \gamma T (P) \right) \Bigr) \leq \delta \Bigl( \partial \gamma(P), \partial L_{n+1} \Bigr) .
\end{align}
Since~${\mathscr{T}}$ is the minimum of the right hand side over all~${n \geq 0}$ and~${\gamma \in T^n}$, and the left hand side is 
equal to~${\delta \Bigl( \partial P, \partial \left( T (P) \right) \Bigr)}$, which is~${\geq \mathscr{T}}$ by definition of~${\mathscr{T}}$, 
we conclude that 
\begin{align}
    \nonumber & \mathscr{T} = \delta \Bigl( \partial P, \partial \left( T (P) \right) \Bigr) .
\end{align}

\begin{figure}[t!]
    \begin{center}
        \includegraphics[width=0.35\textwidth]{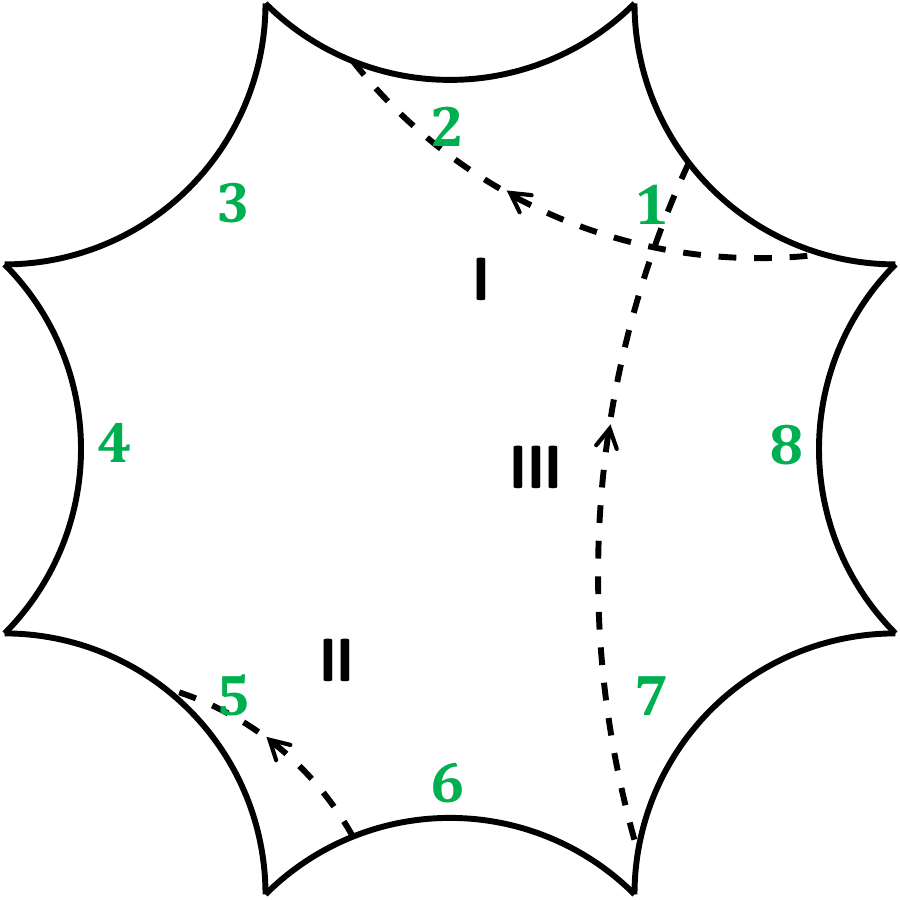}
        \caption{Type I, II, and III segments (dotted black lines), and edge labels (green).}
        \label{F7}
    \end{center}
\end{figure}

Next, we introduce edge labels, which are positive integers modulo the number~${N}$ 
of distinct (in the quotient sense) sides of~${L_0}$. First, the label~${1}$ is assigned 
to an edge of~${L_0}$ without loss of generality. Then, the~${N-1}$ edges to the left of 
this edge are labeled~${2,3, \dots, N-1}$ (Fig.~\ref{F7}). Finally, these labels are copied onto all 
images of~${L_0}$ via the elements of the Fuchsian group~${\mathscr{F}}$. Thus, edges are distinct in 
the quotient sense if and only if they have different labels. The process of filling in labels, which 
we call \emph{label chasing}, will be useful later. 

Now we will prove the main result of this section: 
\begin{align} \label{ET}
    & \mathscr{T} \geq \min \left( \big\{ \delta_{ij} ,~ |i-j| \neq 0,1 \big\} \bigcup \big\{ 2 \epsilon_k \big\} \right) ,
\end{align}
where~${\delta_{ij}}$ denotes the distance between edges~${i}$ and~${j}$ which are distinct and 
nonadjacent, and~${\epsilon_k}$ denotes the 
minimal distance from the midpoint of edge~${k}$ to either of the adjacent edges~${k \pm 1}$. 

\begin{figure}[t!]
    \begin{center}
        \includegraphics[width=0.4\textwidth]{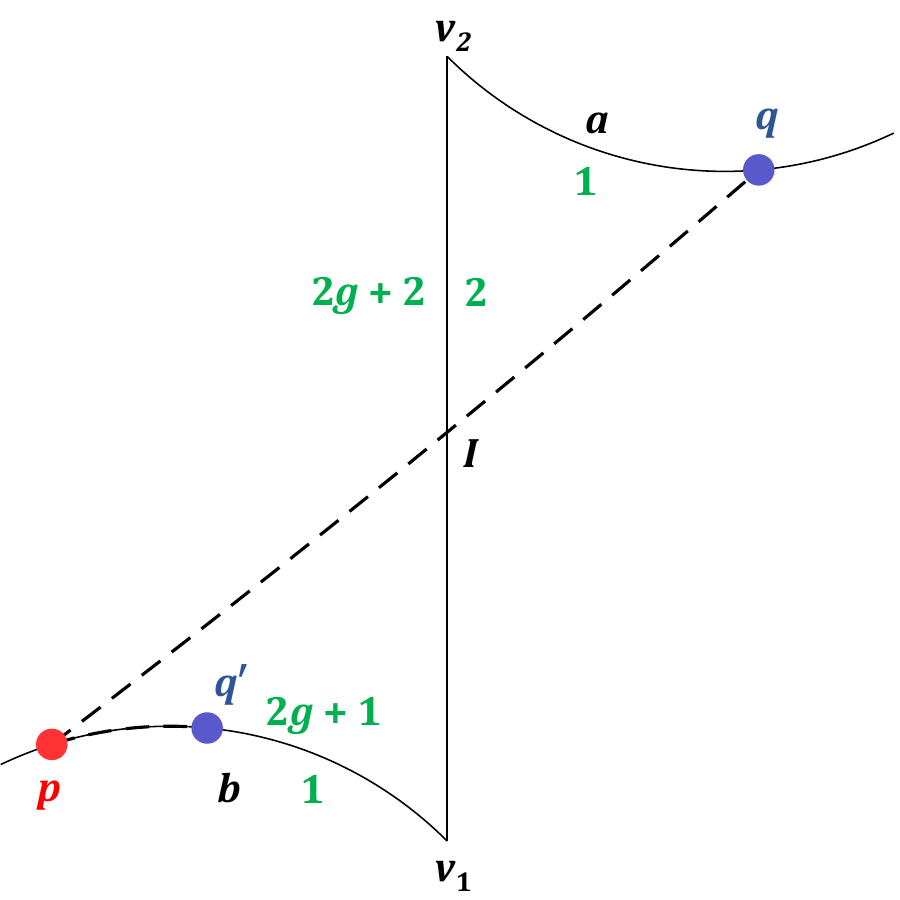}
        \caption{Points~${p,q}$ and image~${q^\prime}$ in case (1) considered in Appendix~\hyperref[Sformula_Bolza]{B}. 
        In this case, the geodesic~${\overline{pq}}$ crosses from one polygon into another through the 
        edge~${\overline{v_1v_2}}$. The intersection point is denoted by~${I}$. 
        The figure shows the relevant edges from both polygons, including 
        their common edge~${\overline{v_1v_2}}$. Segment~${\overline{pI}}$ 
        is a type~I segment while~${\overline{qI}}$ is a type~II segment. 
        Distances~${a = \delta(q,v_2)}$ and~${b = \delta(p,v_1)}$ and edge labels are also shown.}
        \label{FAB1}
    \end{center}
\end{figure}

Starting once more with arbitrary points~${p \in \partial L_0}$ and~${q \in \partial L_1}$, we 
break~${\overline{pq}}$ into smaller segments so that each segment is contained in a polygon~${\gamma(P)}$ 
and has both endpoints on its boundary. For such a segment, suppose the endpoints lie on edges with 
labels~${i,j}$. We categorize a segment (see Fig.~\ref{F7}) as 
\begin{enumerate}
    \item \emph{type I} if~${i-j = 1}$, 
    \item \emph{type II} if~${j-i = 1}$, and 
    \item \emph{type III} if~${|i-j| > 1}$. 
\end{enumerate}
The segments cannot be all type~I, since this would imply that~${\overline{pq}}$ circles around a 
vertex indefinitely without ever reaching~${\partial L_1}$, which cannot happen. Similarly, 
the segments cannot all be type~II. This leaves two cases: there is either (1) a 
consecutive pair of type~I and type~II segments in either order, or (2) a type~III segment. 
In the first case, we can show using elementary geometry that 
a lower bound on the combined length of the type~I and type~II segments, yielding a lower bound 
on~${\mathscr{T}}$, is 
\begin{align}
    & \mathscr{T} \geq \min \big\{ 2 \epsilon_k \big\} .
\end{align}
In the second case,~${\mathscr{T}}$ is lower bounded by the minimum length of a type~III segment, 
\begin{align}
    & \mathscr{T} \geq \min \big\{ \delta_{ij} ,~ |i-j| \neq 0,1 \big\} .
\end{align}
Therefore, in either case, Eq.~\eqref{ET} holds. 

\section{$k_{\min}=1$ for the Bolza surface} \label{Sformula_Bolza}

In this section, we prove that~${k_{\min} = 1}$ for the Bolza surface. 

For an arbitrary point~${z \in P}$, it suffices to show that for any point~${w}$ not in the 
first shell~${L_1}$, defined in Sec.~\ref{Sdistance_formula}, there 
exists a group element~${\gamma}$ such that~${\delta(z,\gamma w) < \delta(z,w)}$. Equivalently, this means 
that the Dirichlet polygon of~${z}$ is contained inside~${L_1}$. 

We prove this by contradiction. 
Assume that there exist points~${z \in P, w \notin L_1}$ that contradict the above 
assertion:~${\delta(z,w) \leq \delta(z,\gamma w)}$ for all~${\gamma}$. This means 
that~${w}$ is either inside the Dirichlet polygon of~${z}$ or on its boundary. 
Using the argument from Appendix~\hyperref[AT_bound]{A}, the 
geodesic~${\overline{zw}}$ must contain either 
\begin{itemize}[leftmargin=1.8cm]
    \item[Case (1):] a segment composed of a pair of type~I and type~II segments in either order, or 
    \item[Case (2):] a type~III segment. 
\end{itemize}
Denote the endpoints of this segment~${p,q}$ such that~${z}$ is closer to~${p}$ 
than to~${q}$. Since the Dirichlet polygon of~${z}$ is convex, contains~${z}$ and~${w}$, 
and~${z,p,q,w}$ are on a line,~${q}$ must be 
in the interior of the Dirichlet polygon of~${z}$. This implies that~${z}$ is in the interior 
of the Dirichlet polygon of~${q}$, and so by the same argument~${p}$ must lie in the interior of 
the Dirichlet polygon of~${q}$. This means that 
\begin{align}\label{Epq}
    & \delta(p,q) < \delta(p, \gamma q) ~\text{for}~ \forall \gamma \neq 1 .
\end{align}

\begin{figure}[t!]
    \begin{center}
        \includegraphics[width=0.4\textwidth]{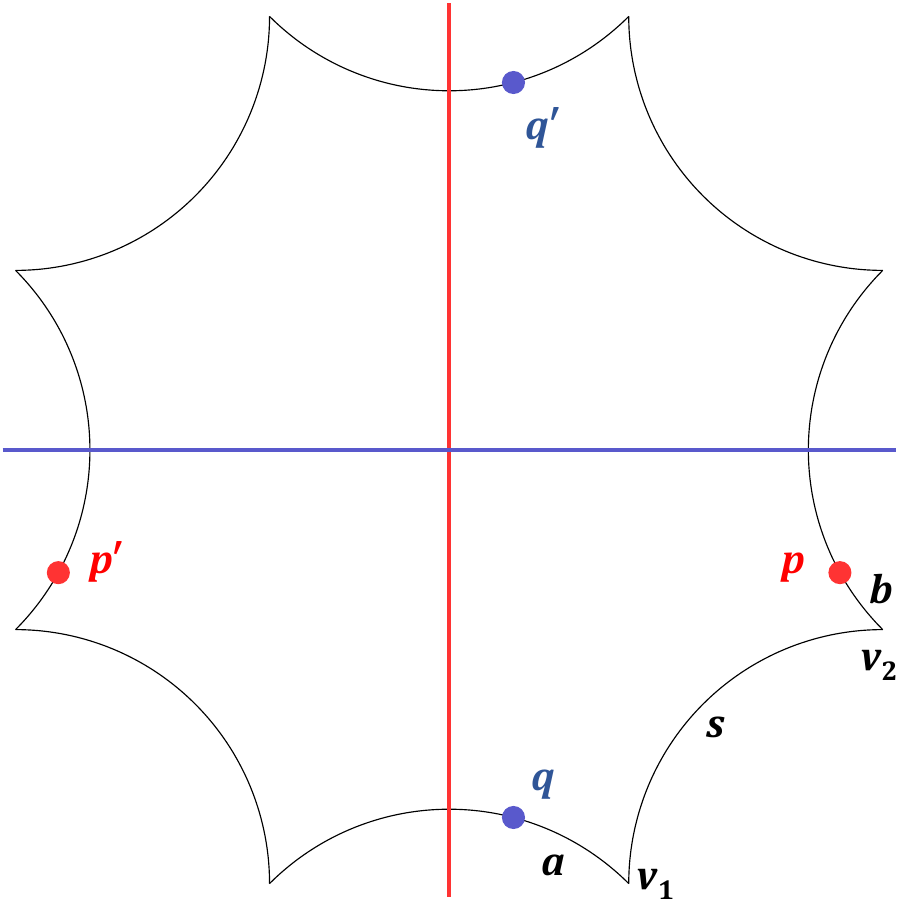}
        \caption{Points~${p,q}$ and images~${p^\prime,q^\prime}$ in case (2). In this case,~${p}$ 
        and~${q}$ are on the same polygon. Perpendicular 
        bisectors of~${\overline{pp^\prime}}$ (red) and~${\overline{qq^\prime}}$ (blue), vertices~${v_1}$ 
        and~${v_2}$ nearest to~${p}$ and~${q}$, distances~${a = \delta(q,v_1)}$ and~${b = \delta(p,v_2)}$, 
        and side length~${s}$ are also shown.}
        \label{FAB2}
    \end{center}
\end{figure}

Consider now case~(1) above, where~${p}$ and~${q}$ are endpoints of a type~I segment joined to a type~II 
segment (Fig.~\ref{FAB1}). Assume without loss of generality that~${q}$ is on side~${1}$. Let~${a = \delta(q,v_2)}$ 
and~${b = \delta(p,v_1)}$. Label-chasing, we find that there must be an image~${q^\prime}$ of~${q}$ 
located on the same side as~${p}$. Applying the triangle 
inequality to triangles~${\triangle v_1Ip}$ and~${\triangle v_2Iq}$, where~${I}$ 
is the intersection between~${\overline{pq}}$ and~${\overline{v_1v_2}}$, gives 
\begin{align}
    \begin{split}
        \delta(p,q) 
        & = \delta (p,I) + \delta (q,I) \\
        & \geq |\delta(v_1,I) - b| + |\delta(v_2,I) - a| \\
        & \geq |\delta(v_1,I) - b + \delta(v_2,I) - a| \\
        & = |s - a - b| \\
        & = \delta(p,q^\prime) ,
    \end{split}
\end{align}
which is a contradiction with Eq.~\eqref{Epq}. 

In case~(2), points ${p}$ and~${q}$ are endpoints of a type~III segment, so they lie on sides of the same 
polygon. Consider the images~${p^\prime,q^\prime}$ 
of~${p,q}$ located on the opposite sides of the polygon, and draw the perpendicular bisectors 
of~${\overline{pp^\prime}}$ and~${\overline{qq^\prime}}$ (Fig.~\ref{FAB2}). 
Observe that~${q}$ must be closer to~${p}$ than to~${p^\prime}$ (on the same side of the vertical perpendicular 
bisector as~${p}$), and~${p}$ must be closer to~${q}$ than to~${q^\prime}$ (on the same side of the horizontal 
perpendicular bisector as~${q}$). This is only possible when~${p}$ and~${q}$ are exactly 2 edges apart and lie in 
the same quarter-polygon, as shown in Fig.~\ref{FAB2}. 

\begin{figure}[t!]
    \begin{center}
        \includegraphics[width=0.4\textwidth]{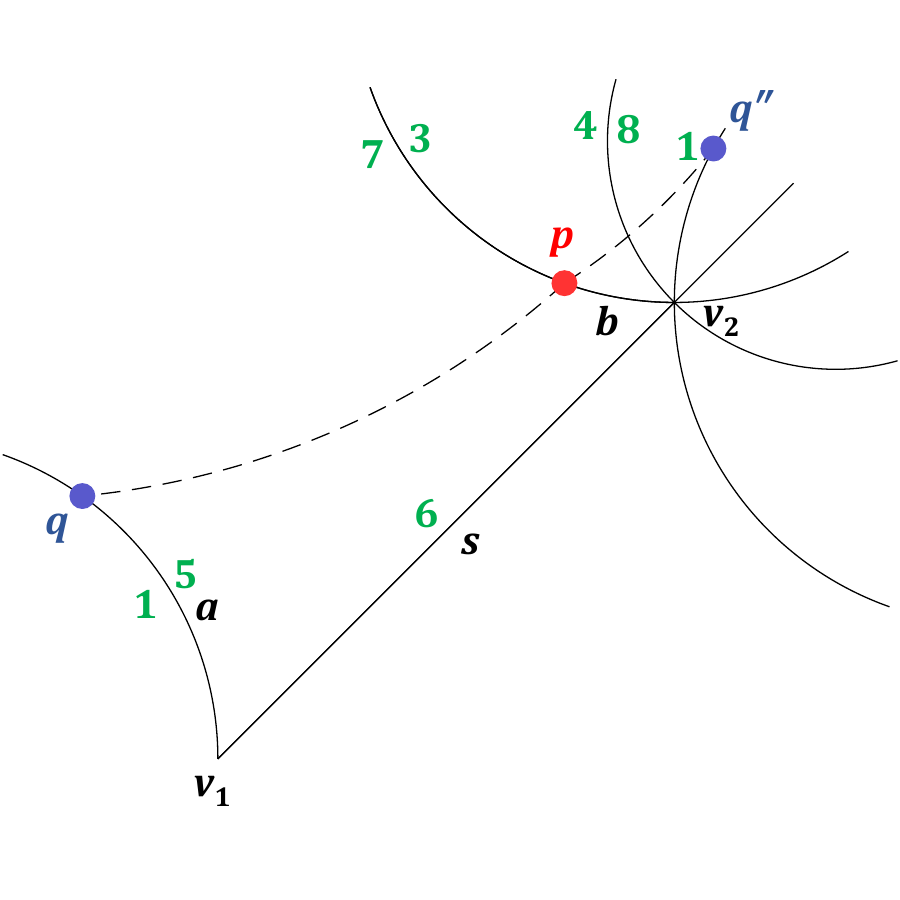}
        \caption{Points~${p,q}$ and image~${q^{\prime\prime}}$ in case (2). Vertices~${v_1,v_2}$, 
        distances~${a,b}$, side length~${s}$, and edge labels are also shown.}
        \label{FAB3}
    \end{center}
\end{figure}

Next, assume without loss of generality that~${q}$ is on side~${1}$ (Fig.~\ref{FAB3}), and use 
label-chasing to determine the location of another image of~${q}$, denoted~${q^{\prime\prime}}$. 
Let~${a = \delta(q,v_1)}$ and~${b = \delta(p,v_2)}$. From Fig.~\ref{FAB2}, we have~${0 < \{a,b\} < s/2}$. 
Applying hyperbolic trigonometry to quadrilateral~${pqv_1v_2}$ and right 
triangle~${\triangle pv_2q^{\prime\prime}}$, one can show that 
\begin{align}\label{Eab}
    \begin{split}
        \cosh (\delta(p,q)) & = \cosh^2 \left( \frac{s}{2} \right) \cosh \left( \frac{s}{2} - b \right) \\
        & \times \cosh \left( \frac{s}{2} - a \right) \\
        & - \sinh \left( \frac{s}{2} \right) \sinh (s-a-b) ,
    \end{split}\\
    \label{Etemp1}
    \cosh (\delta(p,q^{\prime\prime})) & = \cosh a \cosh b .
\end{align}
Rearranging the first equation gives 
\begin{align}
    \nonumber \frac{\cosh (\delta(p,q))}{\cosh b} = & \sinh \left( \frac{s}{2} \right) \cosh \left( \frac{s}{2} - a \right) \\
    \nonumber & \times \bigg[ \sinh \left( \frac{s}{2} \right) \tanh \left( \frac{s}{2} - a \right) \\
    \nonumber & + \cosh \left( \frac{s}{2} \right) - \cosh^2 \left( \frac{s}{2} \right) \bigg] \tanh b \\
    & + C ,
\end{align}

\noindent
where~${C}$ is a function of~${s}$ and~${a}$ only. The coefficient in front of~${\tanh b}$ is a decreasing function 
of~${a}$ and is negative at~${a=0}$. It is therefore negative for 
all~${0 \leq a \leq s/2}$. It follows that as a function of~${b}$,~${\cosh (\delta(p,q)) / \cosh b}$ 
is minimized at~${b = s/2}$. It also follows from Eq.~\eqref{Etemp1} that ${\cosh (\delta(p,q^{\prime\prime})) / \cosh b}$ does not 
depend on~${b}$. To obtain a contradiction, it therefore suffices to show that 
\begin{align}
    & \frac{\cosh (\delta(p,q))}{\cosh b} \geq \frac{\cosh (\delta(p,q^{\prime\prime}))}{\cosh b}
\end{align}
for~${0 \leq a \leq s/2}$ and~${b = s/2}$. Plugging~${b = s/2}$ into Eq.~\eqref{Eab} gives 
\begin{align}
    \begin{split}
        \frac{\cosh (\delta(p,q))}{\cosh (s/2)} & = \cosh \left( \frac{s}{2} \right) \cosh \left( \frac{s}{2} - a \right) \\
        & - \tanh \left( \frac{s}{2} \right) \sinh \left( \frac{s}{2} - a \right) .
    \end{split}
\end{align}
Using~${\cosh (s/2) > 1}$ in the last equation, we get 
\begin{align}
    \begin{split}
        \frac{\cosh (\delta(p,q))}{\cosh (s/2)} & \geq \cosh \left( \frac{s}{2} \right) \cosh \left( \frac{s}{2} - a \right) \\
        & - \sinh \left( \frac{s}{2} \right) \sinh \left( \frac{s}{2} - a \right) \\
        & = \cosh a \\
        & = \frac{\cosh (\delta(p,q^{\prime\prime}))}{\cosh (s/2)} , 
    \end{split}
\end{align}
which is a contradiction completing the proof.

\vspace{1cm}

\begin{acknowledgments}
We thank Florent Balacheff, Maxime Fortier Bourque, Antonio Costa, 
Benson Farb, Svetlana Katok, Gabor Lippner, Marissa Loving, Curtis McMullen, 
Hugo Parlier, John Ratcliffe, and Chaitanya Tappu for 
useful discussions and suggestions. This work was supported by NSF grant 
Nos.~IIS-1741355 and~CCF-2311160.
\end{acknowledgments}

\bibliographystyle{prx}
\bibliography{distance}

\end{document}